\documentclass[12pt]{amsart}
\usepackage[T1]{fontenc}
\usepackage{lmodern}
\usepackage[letterpaper,margin=1.08in]{geometry}
\usepackage{microtype}
\usepackage{mathtools,amssymb}
\usepackage{mathrsfs}
\usepackage{xcolor}
\definecolor{linkblue}{RGB}{0,0,255}

\usepackage[colorlinks=true,linkcolor=linkblue,citecolor=linkblue,
  urlcolor=linkblue,bookmarksnumbered=true,backref]{hyperref}
\hypersetup{pdftitle={Theta singularities and reduced base loci on abelian varieties},
  pdfsubject={Principally polarized abelian varieties, generic vanishing, and base schemes}}
\numberwithin{equation}{section}
\newtheorem{theorem}{Theorem}[section]
\newtheorem{maintheorem}[theorem]{Theorem}
\newtheorem{proposition}[theorem]{Proposition}
\newtheorem{lemma}[theorem]{Lemma}
\newtheorem{corollary}[theorem]{Corollary}
\theoremstyle{definition}

\newtheorem{example}[theorem]{Example}
\theoremstyle{remark}
\newtheorem{remark}[theorem]{Remark}
\newcommand{\OO}{\mathcal O}
\newcommand{\CC}{\mathbf C}
\newcommand{\RR}{\mathbf R}
\newcommand{\LL}{\mathbf L}
\newcommand{\PP}{\mathbf P}
\newcommand{\Pic}{\operatorname{Pic}}
\newcommand{\Sing}{\operatorname{Sing}}
\newcommand{\codim}{\operatorname{codim}}
\newcommand{\ch}{\operatorname{ch}}

\newcommand{\Ann}{\operatorname{Ann}}
\newcommand{\Supp}{\operatorname{Supp}}
\newcommand{\IT}{\mathrm{IT}_0}
\newcommand{\GV}{\mathrm{GV}}
\newcommand{\Db}{D^b_{\mathrm{coh}}}

\title[GV-subschemes and hyperelliptic Jacobians]{GV-subschemes and hyperelliptic Jacobians}
\author{Songsong Huang}
\author{Zhi Jiang}

\address{\rm Shanghai Center for Mathematical Sciences, Fudan University, Shanghai 200438, China}
\email{sshuang22@m.fudan.edu.cn}

\address{\rm Shanghai Center for Mathematical Sciences, Fudan University, Shanghai 200438, China}
\email{zhijiang@fudan.edu.cn} 
\date{}

\begin{document}
\begin{abstract}
We provide a solution to the hyperelliptic Schottky problem: a principally polarized complex abelian variety whose theta
divisor has singular locus of codimension three is a hyperelliptic
Jacobian. As a key step towards this criterion for hyperelliptic Jacobians, we prove a conjecture of Pareschi and Popa which predicts that if the ideal sheaf $I_X$ of a subscheme $X$ of an indecomposable PPAV $(A, \Theta)$ satisfies that $I_X(2\Theta)$ is $\IT$, then $X$ is a GV-subscheme. 

As a corollary of the main result, the base scheme of an ample line bundle on a complex abelian variety is generically reduced along every irreducible component of codimension $2$. This corollary completes the algebro-geometry proof by Arbarello-Codogni-Pareschi in \cite{ACP21, ACP25} of Novikov's conjecture, which characterizes Jacobians via the KP equation.
\end{abstract}
\maketitle
\setcounter{tocdepth}{1}
\tableofcontents

\section{Introduction}
We work over complex number field $\mathbf C$.

The Schottky problem asks how to recognize Jacobians among principally
polarized abelian varieties.  Since the work of Andreotti--Mayer
\cite{AM67}, the dimension of the singular locus of the theta divisor
has been a central invariant in this subject.  Its relation to the
geometry of curves and Prym varieties was developed further by
Beauville \cite{Beauville77} and Debarre \cite{Debarre88}. 
For a Jacobian of a smooth projective curve $C$, Riemann's singularity theorem relates theta
singularities with special divisors of $C$; in particular, a hyperelliptic
Jacobian of dimension $g\geq3$ has its singular locus of dimension
$g-3$.  On the other hand, Ein--Lazarsfeld \cite{EL97} proved that
an irreducible theta divisor is normal and has rational singularities,
so its singular locus has codimension at least three in the ambient
abelian variety.  In particular, the condition
$\codim_A\Sing(\Theta)=2$ characterizes nontrivial products of
principally polarized abelian varieties.  Our first result characterizes
the case $\codim_A\Sing(\Theta)=3$, as formulated in
Debarre's note \cite{Debarre06}.  Beauville kindly explained to us that
he had initially formulated this conjecture in a draft of his thesis and mentioned it in many talks at the time.

\begin{maintheorem}\label{thm:main}
Let $(A,[\Theta])$ be a principally polarized complex abelian variety.  If
\[
 \codim_A\Sing(\Theta)=3,
\]
then there is a smooth hyperelliptic curve $C$ and an isomorphism of
principally polarized abelian varieties
\[
 (A,[\Theta])\simeq(J(C),[\Theta_C]).
\]
\end{maintheorem}

The hypothesis of the above Theorem forces
$\dim A\geq3$.  Beauville \cite{Beauville77} proved this characterization for abelian varieties of
 dimensions at most five.  We 
treat the  higher dimension cases combining generic vanishing, the higher
multiplier ideals of Schnell--Yang \cite{SY26}, and Chen's criterion
for hyperelliptic Jacobians \cite{Chen26}.

For an ample line bundle $L$ on an abelian variety $A$, its base scheme
$\operatorname{Bs}|L|$ is defined by the ideal
\[
 \mathfrak b_L=\operatorname{im}\bigl(H^0(A,L)\otimes L^{-1}
                         \longrightarrow\OO_A\bigr).
\]

\begin{maintheorem}\label{thm:base}
Let $L$ be an ample line bundle on a complex abelian variety $A$.  The
base scheme $\operatorname{Bs}|L|$ is generically reduced along every
irreducible component of codimension two in $A$.
\end{maintheorem}

Debarre showed in \cite{Debarre06} that Theorem~\ref{thm:base} is a consequence of
Theorem~\ref{thm:main}. The codimension bound is sharp:
Example~\ref{ex:nonreduced-base} gives an ample line bundle on an abelian
threefold whose base scheme consists of nine nonreduced double points.
The reducedness assertion  of base loci was discussed in \cite{ACP21} and \cite{ACP25}.
Together with the arguments of Arbarello--Codogni--Pareschi
\cite{ACP21,ACP25}, Theorem~\ref{thm:base} completes an
algebro-geometric proof of Novikov's conjecture, which characterizes
Jacobians via the KP equation.

The main cohomological input is a criterion relating vanishing after
twisting an ideal sheaf by $2\Theta$ to generic vanishing after twisting
by $\Theta$.  Meng \cite{Meng26} proved this implication for ideal
sheaves of effective divisors.  The following theorem extends that
implication to arbitrary nonzero proper coherent ideal sheaves on an
indecomposable principally polarized abelian variety.

\begin{maintheorem}\label{thm:bridge}
Let $(A,[\Theta])$ be an indecomposable principally polarized complex
abelian variety.  Let
$0\ne J\subsetneq\OO_A$ be a coherent ideal sheaf.  If $J(2\Theta)$
satisfies $\IT$, then $X=V(J)$ is a geometrically nondegenerate $\GV$
subscheme.  In particular, $X$ is reduced, irreducible, and
Cohen--Macaulay.  If $g=\dim A$, $d=\dim X$, and
$\theta=c_1(\OO_A(\Theta))$, then
\[
 \chi(\OO_X(\Theta))=1,
 \qquad [X]=\frac{\theta^{g-d}}{(g-d)!}.
\]
\end{maintheorem}

This criterion also bears on the study of subvarieties of minimal
cohomology class.  The curve case is governed by the Matsusaka--Ran
criterion, and Ran established the classification in ambient dimension
four \cite{Ran81}.  Debarre \cite{Debarre95} classified such subvarieties
in Jacobians and formulated a general classification conjecture.
H\"oring \cite{Horing10} established the corresponding uniqueness
statement for minimal surfaces in the intermediate Jacobian of a
generic cubic threefold.  Pareschi--Popa \cite{PP08} connected this
geometric problem with generic vanishing.   Conjecture~A in \cite{PP08} concerns
an indecomposable principally polarized abelian variety $(A,[\Theta])$
of dimension $g$ and a nonempty geometrically nondegenerate reduced closed
subscheme $X\subset A$ of pure dimension $1\leq d\leq g-2$.  With
$\theta=c_1(\OO_A(\Theta))$, it predicts the equivalence of the
following conditions:
\begin{enumerate}
\renewcommand{\labelenumi}{(\arabic{enumi})}
\item $[X]=\theta^{g-d}/(g-d)!$;
\item $I_X(\Theta)$ is a $\GV$ sheaf;
\item $I_X(2\Theta)$ satisfies $\IT$;
\item $\OO_X(\Theta)$ is $M$-regular and
      $\chi(\OO_X(\Theta))=1$.
\end{enumerate}
Here $M$-regularity means that $\codim V^i(\OO_X(\Theta))>i$ for
all $i>0$; the other cohomological conditions are recalled in
Subsection~\ref{subsec:notation}.  The fifth condition asserts that
$(A,[\Theta])$ is a polarized Jacobian of a smooth projective curve $C$ and $X$ is a translate of
$\pm W_d(C)$, or that $g=5$, $d=2$, and $(A,[\Theta])$ is the
intermediate Jacobian of a smooth cubic threefold $Y$, with $X$ a translate
of $\pm F(Y)$, where $F(Y)$ is the Fano variety of lines of $Y$.  These examples
satisfy all four conditions \cite{PP08,Horing07}.

Pareschi--Popa \cite{PP08} proved that (2) implies (1)
 and obtained the geometric classification for curves and
subschemes of codimension two.  Lombardi--Tirabassi \cite{LT15}
showed that the inclusion of a generating reduced $\GV$ subscheme
cannot factor through a nontrivial isogeny.
Schreieder \cite{Schreieder16,Schreieder17}
treated $\GV$ subvarieties with a curve summand and established
structural results for $\GV$ subschemes and their theta-duals.
Casalaina-Martin--Popa--Schreieder \cite{CMPS18} proved the classification
of geometrically nondegenerate $\GV$ subschemes of a PPAV of dimension $5$.
Jiang--Pareschi \cite{JP20} determined the Hilbert polynomial and
structure-sheaf cohomology of geometrically nondegenerate $\GV$
subvarieties.  They also proved that geometrically nondegenerate
$\GV$ surfaces in indecomposable principally polarized abelian
varieties are normal with rational singularities.
More recently, Chen \cite{Chen26} classified geometrically nondegenerate
$\GV$ subschemes of codimension three and, by theta-duality,
geometrically nondegenerate $\GV$ surfaces, in ambient dimension at
least six.  Meng \cite{Meng26} proved the divisor analogue on an
indecomposable principally polarized abelian variety: for a nonzero
effective divisor, the three cohomological conditions are equivalent
to being a translate of the theta divisor.  Theorem~\ref{thm:bridge}
gives the following consequence
in the setting of Conjecture~A.

\begin{corollary}\label{cor:conjecture-a}
Let $X\subset A$ be a subscheme. Conditions~\textup{(2)},
\textup{(3)}, and~\textup{(4)} are equivalent, and each implies
condition~\textup{(1)} and that $X$ is an integral and geometric non-degenerate subvariety of $A$.
\end{corollary}

\begin{proof}
Pareschi--Popa \cite{PP08} proved
$(2)\Longleftrightarrow(4)$ and $(2)\Longrightarrow(3)$.   The implication
$(3)\Longrightarrow(2)$ follows from Theorem~\ref{thm:bridge} applied
to $J=I_X$.  The rest  also follows from Theorem~\ref{thm:bridge}. 
\end{proof}

Thus the theorem establishes the equivalence of the cohomological
conditions in Conjecture~A.  The reverse implication from minimal
class to these conditions is not supplied by our argument.

We briefly describe the proof of Theorem~\ref{thm:bridge}.
The Fourier--Mukai transform gives
$V^i(J(\Theta))\subseteq V^0(J(\Theta))$ for all $i$.
Using addition of subvarieties and an induction on the dimension, we show that
$\codim V^0(J(\Theta))\geq d+1$.  Together with
$V^i(J(\Theta))=\varnothing$ for $i>d+1$, this proves that $J(\Theta)$ is a GV sheaf.  The Chern character calculation in
Lemma~\ref{lem:chern} then shows that $X$ is  geometric nondegenerate, and
\cite{PP08,Schreieder17} supplies the remaining assertions.

To prove Theorem~\ref{thm:main}, we apply this criterion to a higher
multiplier ideal from \cite{SY26}.  Its zero scheme is a
codimension-three $\GV$ subscheme contained in $\Sing(\Theta)$.
Theta-duality then produces a singular surface, to which Chen's
criterion \cite{Chen26} applies.  We recall the necessary background
in Section~\ref{sec:preliminaries} and prove Theorem~\ref{thm:bridge}
in Sections~\ref{sec:retraction} and~\ref{sec:codimension}. We prove Theorem~\ref{thm:main}
in Section~\ref{sec:theta} and
show that the codimension bound in Theorem~\ref{thm:base} is sharp
in Section~\ref{sec:base-loci}.

\subsection*{Acknowledgements}
The authors thank Olivier Debarre and Arnaud Beauville for helpful
discussions and valuable suggestions and Giuseppe Pareschi for many
useful comments. Mihnea Popa
kindly informed us that he and Fanjun Meng have recently obtained similar results independently.

The authors are supported by the Foundation for Innovative Research Groups
of the National Natural Science Foundation of China (No.~12121001).
The second author is a member of the Key Laboratory of Mathematics for
Nonlinear Science, Fudan University, and is also supported by the National
Natural Science Foundation of China (No.~12371044), the Shanghai Leading
Talent Program of Eastern Talent Plan, and the Natural Science Foundation
of Shanghai (No.~23ZR1404000).
Besides assisting with editing and revising the manuscript, AI tools suggested us to complete the proof of Theorem \ref{thm:main} by applying Theorem \ref{thm:chen}   and  constructed the counterexample in Example~\ref{ex:nonreduced-base}.
The authors take full responsibility for the mathematical content and
the final text.

\section{Preliminaries}
\label{sec:preliminaries}

We work over $\CC$.  We recall the definitions and results used in the
proofs below.  For an abelian variety $A$, all tangent spaces are
identified by translation, and $\widehat A$ denotes $\Pic^0(A)$.

\subsection{Notation and generic vanishing}\label{subsec:notation}
We first fix the cohomological and geometric nondegeneracy conditions
that enter our generic vanishing criterion.
For a coherent sheaf $F$ on an abelian variety $A$, let
\[
 V^i(F)=\{P\in\Pic^0(A):H^i(A,F\otimes P)\ne0\}.
\]
The sheaf $F$ satisfies $\IT$ if $V^i(F)=\varnothing$ for every $i>0$,
and it is a $\GV$ sheaf if $\codim V^i(F)\geq i$ for every $i>0$.
A closed subscheme $X\subset A$ is a $\GV$ subscheme of $(A,[\Theta])$
if $I_X(\Theta)$ is $\GV$.
These conditions are preserved under tensoring with a line bundle in
$\Pic^0(A)$, since such a twist translates the cohomological support loci in $\widehat{A}$.

If $d=\dim X$, geometric nondegeneracy means that no nonzero decomposable
holomorphic $d$-form on $A$ pulls back to zero on resolutions of all the
$d$-dimensional components of $X_{\mathrm{red}}$.  We use this convention
also when $X$ is initially nonreduced or nonpure; see
\cite[\S2.4]{Schreieder17}.  Tangent spaces to subvarieties of $A$
are identified with subspaces of $T_0A$ by translation.  Chern character
identities are taken in rational cohomolog.  Tensor products and
pullbacks of derived objects are derived; tensoring by a line bundle and
external products of coherent sheaves over $\CC$ need no additional
derived notation.

\subsection{Addition of subvarieties and theta divisors}

We recall some results on the addition of subvarieties and theta divisors on abelian varieties.

\begin{theorem}
\label{thm:addition}
Let $V,W$ be irreducible subvarieties of an abelian variety $A$ of
 dimension $g$.
\begin{enumerate}
\renewcommand{\labelenumi}{(\roman{enumi})}
\item Let $K$ be the identity component of the stabilizer of $V+W$,
and let $p:A\to A/K$ be the quotient.  Then
\[
 \dim p(V+W)=\dim p(V)+\dim p(W).
\]
In particular, if $V+W$ has finite stabilizer, then
$\dim(V+W)=\dim V+\dim W$.
\item If $V$ is geometrically nondegenerate, then
\[
 \dim(V+W)=\min\{\dim V+\dim W,g\}.
\]
\end{enumerate}
\end{theorem}

\begin{proof}
Part~(i) is \cite[Chapitre~VIII, théorème~2.2]{Debarre99}.
Part~(ii) is \cite[Corollaire~2.7]{Debarre95Connectivity}.
\end{proof}

\begin{theorem}
\label{thm:theta-normal}
The theta divisor of a principally polarized complex abelian variety
$(A,[\Theta])$ is irreducible if and only if $(A,[\Theta])$ is
indecomposable.
An irreducible theta divisor is normal and has rational singularities.
Its singular locus, when nonempty, has codimension at least $3$ in $A$.
\end{theorem}

\begin{proof}
The irreducibility assertion follows from the decomposition theorem
in \cite[Theorem~4.3.1]{BL04}.  Normality and rational
singularities are due to Ein--Lazarsfeld \cite[Theorem~1]{EL97}.
For the last assertion, a normal variety is regular in codimension one,
so $\Sing(\Theta)$ has codimension at least two in $\Theta$.
Since $\Theta$ is a divisor in $A$, its singular locus has codimension
at least three in $A$.
\end{proof}

\begin{theorem}
\label{thm:summands}
Let $A$ be an abelian variety of dimension $g$, and let $V,W\subset A$
be irreducible subvarieties of positive dimensions $d$ and $g-d-1$.
If $V+W$ is an ample divisor with rational singularities, then $V$ and
$W$ are geometrically nondegenerate.
\end{theorem}

\begin{proof}
This is Schreieder's Theorem (\cite[Theorem~1(1)]{Schreieder17}).
\end{proof}

\subsection{Theta-dual}

Theta-dual records the translates of a theta divisor containing a
subscheme.  We need its scheme-theoretic formulation and its behavior
on geometrically nondegenerate $\GV$ subschemes.

\begin{proposition}
\label{prop:theta-dual}
Let $(A,[\Theta])$ be a principally polarized abelian variety and let
$X\subset A$ be a closed subscheme.  There is a closed subscheme
$S\subset A$ representing the condition $X+a\subseteq\Theta$.
Precisely, for any scheme $Q$ and morphism $a:Q\to A$, the map $a$
factors through $S$ if and only if the relative translation
\[
 X\times Q\longrightarrow A\times Q,\qquad
 (x,q)\longmapsto(x+a(q),q)
\]
factors through $\Theta\times Q$.
\end{proposition}

\begin{proof}
See Gulbrandsen--Lahoz
\cite[Definition~2.3 and Propositions~2.5--2.6]{GL11} for the
representability and universal property of the theta-dual.
\end{proof}

We call $S$ the \emph{theta-dual} of $X$.  This convention agrees with
\cite{Chen26}; it is the inverse of the theta-dual convention in
\cite{Schreieder17}.  The following properties are unchanged by inversion.

\begin{theorem}
\label{thm:gv-duality}
Let $(A,[\Theta])$ be a principally polarized abelian variety of
 dimension $g$, and let $X\subsetneq A$ be a nonempty geometrically
nondegenerate closed $\GV$ subscheme of dimension $d$.
Let $S$ be its theta-dual, and put $\theta=c_1(\OO_A(\Theta))$.
Then $X$ and $S$ are pure-dimensional Cohen--Macaulay schemes, and
$S$ is a geometrically nondegenerate $\GV$ subscheme of dimension
$g-d-1$.  Their cohomology classes are
\[
 [X]=\frac{\theta^{g-d}}{(g-d)!},
 \qquad [S]=\frac{\theta^{d+1}}{(d+1)!}.
\]
If $(A,[\Theta])$ is indecomposable, both $X$ and $S$ are reduced
and irreducible.
\end{theorem}

\begin{proof}
The duality and minimal class assertions are due to Pareschi and Popa
\cite[Theorems~5.2 and~6.1]{PP08}, in the scheme-theoretic form of
Schreieder's \cite[Theorem~8]{Schreieder17}; the integrality assertion is
\cite[Theorem~2(1)]{Schreieder17}.
\end{proof}

\subsection{Higher multiplier ideals}

Higher multiplier ideals provide the subscheme of the singular locus
to which we will apply our generic vanishing criterion.  We record
the containment, vanishing, and support properties needed for this step.

For an effective divisor $D$ on a smooth complex variety, write
$I_{k,\alpha}(D)$ for the higher multiplier ideals of Schnell and
Yang, where $k\geq0$ and $\alpha\in\mathbf Q$.
For fixed $k$, this family is discrete and right continuous
\cite[\S1.1]{SY26}.  Its strict left limit is therefore
\[
 I_{k,<\alpha}(D)=I_{k,\alpha-\varepsilon}(D)
 \qquad(0<\varepsilon\ll1).
\]

\begin{theorem}
\label{thm:higher-multiplier}
\leavevmode
\begin{enumerate}
\renewcommand{\labelenumi}{(\roman{enumi})}
\item Let $D$ be an effective divisor on a smooth complex variety $M$.
If $Z$ is an irreducible component of
$\operatorname{Sing}(D)$
of codimension three in $M$, then
\[
 I_{1,<-1/2}(D)\subseteq I_Z.
\]
\item Let $D$ be an ample effective divisor on an abelian variety $A$,
and put $L=\OO_A(D)$.  For every integer $k\geq0$, rational number
$\alpha\in[-1,0)$, and $P\in\Pic^0(A)$,
\[
 H^i\bigl(A,L^{k+1}\otimes I_{k,\alpha}(D)\otimes P\bigr)=0
 \qquad(i>0).
\]
\item For an effective divisor $D$ on a smooth complex variety $M$,
$k\geq0$, and rational $-1<\alpha<0$, the ideal
$I_{k,\alpha}(D)$ is the unit ideal where $D$ is smooth or absent.
\end{enumerate}
\end{theorem}

\begin{proof}
Part~(i) is Schnell--Yang's \cite[Theorem~1.3]{SY26} with
$(m,k,r,q)=(2,1,1,1)$, part~(ii) is
\cite[Theorem~1.6(3)]{SY26}, and part~(iii) follows from
\cite[Lemma~3.10]{SY26}.
\end{proof}
In particular, the ideal in part~(i) can be
written as $I_{1,-1/2-\varepsilon}(D)$ for sufficiently small rational
$0<\varepsilon<1/2$, so part~(ii) applies to this same ideal when $D$
is ample on an abelian variety.

\subsection{Jacobian criteria}

The final ingredients characterize hyperelliptic Jacobians.  The first
result covers dimensions three to five; in higher dimensions we use
the second criterion, expressed in terms of a singular theta-dual.

\begin{theorem}
\label{thm:beauville}
Let $(A,[\Theta])$ be a principally polarized abelian variety of
 dimension $g$ with $3\leq g\leq5$.  If
$\dim\Sing(\Theta)=g-3$, then $(A,[\Theta])$ is the principally
polarized Jacobian of a smooth hyperelliptic curve.
\end{theorem}

\begin{proof}
See Beauville \cite[Proposition~6.6(2)]{Beauville77}.
\end{proof}

\begin{theorem}
\label{thm:chen}
Let $(A,[\Theta])$ be an indecomposable principally polarized complex
abelian variety of dimension $g\geq6$.  Let $X\subset A$ be a nonempty
$\GV$ subscheme of dimension $g-3$ whose reduction is geometrically
nondegenerate.
If the scheme-theoretic theta-dual $S$ of $X$ is singular, then there
is a smooth hyperelliptic curve $C$ of genus $g$ and an isomorphism
of principally polarized abelian varieties
\[
 (A,[\Theta])\simeq(J(C),[\Theta_C]).
\]
\end{theorem}

\begin{proof}
This is Chen's \cite[Theorem~4.7]{Chen26}.
\end{proof}

\section{Cohomological support loci and Chern characters}
\label{sec:retraction}

We begin the proof of Theorem~\ref{thm:bridge}. After a translation, we may assume that $\Theta$ is symmetric. We write $\iota=(-1)_A$ and identify $A$ with $\widehat A$ by the principal
polarization.  Let
\[
 \Phi:\Db(A)\longrightarrow\Db(A)
\]
be the Fourier--Mukai transform for the normalized Poincar\'e bundle.
With this normalization $\Phi(\OO_A(\Theta))\simeq \OO_A(-\Theta)$; see \cite{Mukai81}.
Put
\[
 K=\Phi\bigl(J(\Theta)\bigr).
\]

\begin{lemma}\label{lem:retraction}
The object $K\otimes \OO_A(-\Theta)$ is a direct summand of
$K\otimes^{\LL}K$ in $\Db(A)$.
\end{lemma}

\begin{proof}
Consider the difference map and
\[
 \delta:A\times A\longrightarrow A,\quad (u,v)\longmapsto v-u,
\]
and its fiber-preserving involution
$$\sigma: A\times A\to A\times A, \quad(u,v)\longmapsto(-v,-u).$$
Set
\[
 E=J(\Theta)\boxtimes\OO_A(\Theta),
 \qquad \mathcal L=\OO_A(\Theta)\boxtimes\OO_A(\Theta).
\]
Note that $E$ is flat over $A$ via $\delta$, and that for each $a\in A$, its restriction to the fiber $\delta^{-1}(a)$ is isomorphic to $J(2\Theta)\otimes P_a$ for some $P_a\in\Pic^0(A)$.  Since $J(2\Theta)$ is $\IT$, cohomology and base change
therefore give vector bundles
\[
 \mathcal B=\RR\delta_*E=\delta_*E,
 \qquad \mathcal W=\delta_*\mathcal L,
\]
and a sheaf injection $\mathcal B\hookrightarrow\mathcal W$.

Choose a normalized symmetry $\iota^*\OO_A(\Theta)\simeq \OO_A(\Theta)$.  Applying it to both
factors and interchanging the factors lifts $\sigma$ to an action
of $\mathcal L$, and hence of $\mathcal W$.  On the fiber over $0$, $\sigma\vert_{\delta^{-1}(0)}$ is nothing but $\iota$, up to $\delta^{-1}(0)\cong A$ and 
the action on $\mathcal{W}$ is exactly the action by $\iota^*$ on $H^0(A,\OO_A(2\Theta))$.  Every section
of $\OO_A(2\Theta)$ is even
\cite[Corollary~4.6.6 and \S4.8]{BL04}; thus the action on that fiber
is the identity.  The positive and negative eigensubbundles have locally
constant ranks.  Since $A$ is connected, the action on all of
$\mathcal W$ is the identity.

The evaluation map $\delta^*\mathcal W\to\mathcal L$ is $\sigma$-equivariant.
Its restriction to $\delta^*\mathcal B$ lands in $E$, and the identity
action on $\mathcal W$ forces the image to lie also in
$\OO_A(\Theta)\boxtimes\iota^*(J(\Theta))$.  Inside $\mathcal L$ we have
\begin{equation}\label{eq:intersection}
 \begin{aligned}
 &\bigl(J(\Theta)\boxtimes\OO_A(\Theta)\bigr)
   \cap\bigl(\OO_A(\Theta)\boxtimes\iota^*(J(\Theta))\bigr)\\
 &\qquad=J(\Theta)\boxtimes\iota^*(J(\Theta)).
 \end{aligned}
\end{equation}
Thus the evaluation
factors as
\[
 \delta^*\mathcal B\longrightarrow J(\Theta)\boxtimes\iota^*(J(\Theta))
          \longrightarrow E.
\]
Derived adjunction gives morphisms
\begin{equation}\label{eq:geometric-retraction}
 \mathcal B\longrightarrow
 \RR\delta_*(J(\Theta)\boxtimes\iota^*(J(\Theta)))
 \longrightarrow\mathcal B
\end{equation}
whose composite is the identity: its adjoint is the original counit
$\delta^*\mathcal B\to E$.

We have
\[
 \Phi\bigl(\RR\delta_*(U\boxtimes V)\bigr)
 \simeq \iota^*\Phi(U)\otimes^{\LL}\Phi(V).
\]
It follows directly from the following formula for the  Poincar\'e
bundle $\mathcal P$ on $A\times \Pic^0(A)$: considering $\delta\times \mathrm{Id}_{\Pic^0(A)}: (A\times A)\times \Pic^0(A)\rightarrow A\times \Pic^0(A)$, then $(\delta\times \mathrm{Id}_{\Pic^0(A)})^*\mathcal P\simeq p_{13}^*\mathcal P^{-1}\otimes p_{23}^*\mathcal P$.
Applying it to
\eqref{eq:geometric-retraction} yields a retraction of
$\iota^*K\otimes \OO_A(-\Theta)$ from
$\iota^*(K\otimes^{\LL}K)$.  Pullback by $\iota$ gives the assertion.
\end{proof}

%\begin{remark}
%The argument uses the injection $\mathcal B\hookrightarrow\mathcal W$
%as a map of sheaves, without asserting injectivity on every ordinary
%fiber.  It also keeps all higher direct images in the middle term of\eqref{eq:geometric-retraction}.  Both points matter when $J$ has arbitrary scheme structure.
%\end{remark}

\begin{corollary}\label{cor:support}
For every $i\geq0$, one has $V^i(J(\Theta))\subseteq V^0(J(\Theta))$.  Moreover,
$V^0(J(\Theta))$ is a proper closed subset of $\widehat A$, and
\begin{equation}\label{eq:chi}
 \chi(J(\Theta))=0,\qquad \chi(\OO_X(\Theta))=1.
\end{equation}
\end{corollary}

\begin{proof}
For $P\in\widehat A$, derived base change identifies the derived fiber
of $K$ with $C_P=\RR\Gamma(A,J(\Theta)\otimes P)$.  This complex has no
negative cohomology.  If it is nonzero and its least nonzero
cohomology degree is $k$, the least nonzero degree of
$C_P\otimes^{\LL}_{\CC}C_P$ is $2k$.  Lemma~\ref{lem:retraction}
makes $C_P$, up to a one-dimensional tensor factor, a direct summand
of this tensor square.  Thus $k=0$, proving the support inclusions.

Since $J$ is proper, choose a closed point $x\in X$.
Every nonzero section of $J(\Theta)\otimes P$ gives the unique, up to scalar,
section of $\OO_A(\Theta)\otimes P$ vanishing at $x$.  Consequently $V^0(J(\Theta))$ is
contained in a theta translate and is proper.  At a point outside
$V^0(J(\Theta))$ all cohomology of $J(\Theta)\otimes P$ vanishes, so $\chi(J(\Theta))=0$.
The exact sequence
\begin{equation}\label{eq:ideal-sequence}
 0\longrightarrow J(\Theta)\longrightarrow \OO_A(\Theta)\longrightarrow
 \OO_X(\Theta)\longrightarrow0
\end{equation}
and $\chi(\OO_A(\Theta))=1$ give the second equality.
\end{proof}

\begin{lemma}\label{lem:chern}
Let $V=V^0(J(\Theta))$.  If $\codim V\geq r$, then
\begin{equation}\label{eq:chern}
 \ch_{g-j}(\OO_X(\Theta))
     =\frac{\theta^{g-j}}{(g-j)!}\qquad(0\leq j<r).
\end{equation}
\end{lemma}

\begin{proof}
Corollary~\ref{cor:support} and derived base change imply that all
cohomology sheaves of $K$ are supported on $V$.  Hence
$\ch_j(K)=0$ for $j<r$.  The cohomological Fourier transform $\Phi^H$ maps
$H^{2(g-j)}(A,\mathbf Q)$ isomorphically to $H^{2j}(A,\mathbf Q)$
and satisfies
\[
 \Phi^H\left(\frac{\theta^{g-j}}{(g-j)!}\right)
       =(-1)^j\frac{\theta^j}{j!},
\]
see \cite[Theorem 16.5.5]{BL04}.
As the Todd class of $A$ is $1$, $\Phi^H$ commutes with Chern characters (\cite[Proposition 16.3.1]{BL04}).
Thus $\ch_{g-j}(J(\Theta))=0$ for $j<r$.  Taking Chern characters in
\eqref{eq:ideal-sequence} proves \eqref{eq:chern}.
\end{proof}

\section{Codimension estimation and generic vanishing}
\label{sec:codimension}

We retain the notation of Section~\ref{sec:retraction} and put
$d=\dim X$.  Under the identification
$a\mapsto t_a^*\OO_A(\Theta)\otimes\OO_A(-\Theta)$ of $A$ and
$\widehat A$, a point $a$ in $V=V^0(J(\Theta))$ means that
$X+a\subseteq\Theta$ scheme-theoretically.
In particular,
\begin{equation}\label{eq:incidence}
 Y+Z\subseteq\Theta
\end{equation}
for every irreducible component $Y$ of $X_{\mathrm{red}}$ and every
irreducible subvariety $Z\subseteq V$.

\begin{proposition}\label{prop:codimension}
One has $\codim V\geq d+1$.
\end{proposition}

\begin{proof}
If $d=0$, properness of $V$ proves the assertion.  Suppose $d>0$.
Assume to the contrary that $\codim V=1$. Let $Z$ be a divisorial component of $V$ and choose a
positive-dimensional component $Y$ of $X_{\mathrm{red}}$.
By \eqref{eq:incidence}, we have $Y+Z\subseteq\Theta$; the two sides have
dimension $g-1$, and $\Theta$ is irreducible by the decomposition
theorem \cite[Theorem~4.3.1]{BL04}, so $Y+Z=\Theta$.
The stabilizer of the ample divisor $\Theta$ is finite.
Theorem~\ref{thm:addition}(i) would give $\dim Y+(g-1)=g-1$, a contradiction.
Thus $\codim V\geq2$.

Assume inductively that $2\leq r\leq d$ and $\codim V= r$.
Let $Z$ be a component of $V$ of codimension $r$, and let
$Y_1,\ldots,Y_s$ be all the components of $X_{\mathrm{red}}$.
If $Y_i+Z=\Theta$ for some $i$, Theorem~\ref{thm:addition}(i)
gives $\dim Y_i=r-1$.  Both summands are positive-dimensional, and
their dimensions add to $g-1$.  The divisor $\Theta$ has rational
singularities by Theorem~\ref{thm:theta-normal}, so
Theorem~\ref{thm:summands} makes $Z$ geometrically nondegenerate.
For a component $Y_{j}$ of dimension $d$,
Theorem~\ref{thm:addition}(ii) would then give
\[
 \dim(Y_{j}+Z)=\min\{d+g-r,g\}=g,
\]
contrary to \eqref{eq:incidence}.  Consequently
\begin{equation}\label{eq:smaller-sums}
 \dim(Y_i+Z)\leq g-2\qquad\text{for every }i.
\end{equation}

Choose one general smooth point $z\in Z$.  At a general smooth point $y\in Y_i$,
considering the generic rank of the addition map,  \eqref{eq:smaller-sums} implies
\begin{equation}\label{eq:dim-quotient-tengent}
 \dim\bigl((T_yY_i+T_zZ)/T_zZ\bigr)\leq r-2.
\end{equation}
Since $\Ann(T_zZ)\subset H^0(A, \Omega_A^1)$ has dimension $r$, choose a nonzero form
\[
 \eta\in\bigwedge^{r-1}\Ann(T_zZ)
          \subset H^0(A,\Omega_A^{r-1}).
\]
The form $\eta$ vanishes on a resolution of each $Y_i$ by \eqref{eq:dim-quotient-tengent}, and hence
on a resolution of every integral subvariety $D\subseteq X_{\mathrm{red}}$.
Let $F=\OO_X(\Theta)$ be a coherent sheaf supported on $X$.
Hence, we have
\begin{equation}\label{eq:zero-pairing}
 \int_A\ch_{g-r+1}(F)\wedge\eta\wedge\overline\eta
 %=\sum_j a_j\int_{\widetilde D_j}
       %f_j^*\eta\wedge f_j^*\overline\eta
       =0.
\end{equation}
On the other hand, Lemma~\ref{lem:chern} applies with $j=r-1$.
It replaces the Chern character component in \eqref{eq:zero-pairing}
by $\theta^{g-r+1}/(g-r+1)!$.  The K\"ahler form representing
$\theta$ gives the strictly positive norm
\[
 \mathrm i^{(r-1)^2}\int_A
 \frac{\theta^{g-r+1}}{(g-r+1)!}
       \wedge\eta\wedge\overline\eta>0.
\]
This contradicts \eqref{eq:zero-pairing}.  There is no component of
codimension $r$, so $\codim V\geq r+1$.  Induction proves the result.
\end{proof}

\begin{remark}\label{rem:pareschi-popa-bound}
Pareschi--Popa prove in \cite[Lemma~4.4]{PP08} that, for a geometrically
nondegenerate reduced subscheme $Y\subset A$ of pure dimension
$e\leq g-2$, one has
\[
 \dim V^0(I_Y(\Theta))\leq g-e-1.
\]
At this stage, however, $X$ is not known to be geometrically
nondegenerate, reduced, or pure-dimensional.  In fact, its geometric
nondegeneracy is deduced below from Proposition~\ref{prop:codimension}
via Equation \eqref{eq:minimal-class}.  Thus their
lemma cannot replace the preceding argument without circularity.
\end{remark}

\begin{proof}[Proof of Theorem~\ref{thm:bridge}]
Since $\OO_A(\Theta)\otimes P$ is ample, \eqref{eq:ideal-sequence} gives
\[
 H^i(A,J(\Theta)\otimes P)=0\qquad(i>d+1)
\]
for every $P\in\Pic^0(A)$.  For $1\leq i\leq d+1$,
Corollary~\ref{cor:support} and Proposition~\ref{prop:codimension}
give
\[
 \codim V^i(J(\Theta))\geq\codim V\geq d+1\geq i.
\]
Hence $J(\Theta)$ is $\GV$.  Lemma~\ref{lem:chern}, with $r=d+1$ and $j=d$,
also gives the top-dimensional cycle
\begin{equation}\label{eq:minimal-class}
 [X]=\ch_{g-d}(\OO_X(\Theta))
          =\frac{\theta^{g-d}}{(g-d)!}.
\end{equation}
If a nonzero decomposable holomorphic $d$-form vanished on resolutions
of all the top-dimensional components of $X_{\mathrm{red}}$, its
pairing with the cycle in \eqref{eq:minimal-class} and its conjugate
would be zero.  The same positive norm argument makes this impossible.
Thus $X$ is geometrically nondegenerate.
Theorem~\ref{thm:gv-duality} now implies that $X$ is reduced,
irreducible, and Cohen--Macaulay.  Equation~\eqref{eq:chi}
proves the remaining assertion.
\end{proof}

\section{The hyperelliptic Jacobian characterization}
\label{sec:theta}

The following observation provides the singularity hypothesis in
Theorem~\ref{thm:chen}.

\begin{lemma}\label{lem:singular-dual}
Let $X\subseteq\Sing(\Theta)$ be a reduced closed subscheme of $A$,
and let $S$ be its theta-dual.  Then $0\in S$ and
\[
 T_0S=T_0A.
\]
In particular, if $\dim_0S<\dim A$, the scheme $S$ is singular at $0$.
\end{lemma}

\begin{proof}
See also the argument in
\cite[proof of Lemma~4.1]{Chen26}.  Since $X\subseteq\Theta$, we have
$0\in S$.  For $v\in T_0A$, let $\partial_v$ be the translation-invariant
derivation of $\OO_A$ whose value at $x\in A$ is
$(dt_x)_0(v)\in T_xA$, where $t_x$ denotes translation by $x$.
Write $a_v:D=\operatorname{Spec}(\CC[\epsilon]/(\epsilon^2))\to A$
for the tangent vector $v$, and consider
\[
 f=m\circ(i\times a_v):X\times D\longrightarrow A,
\]
where $i:X\hookrightarrow A$ is the inclusion and $m$ is addition.
For a local regular function $g$, its pullback is
\[
 f^*(g)=g|_X+\epsilon(\partial_vg)|_X.
\]
For a local equation $h$ of $\Theta$, both $h|_X$ and
$(\partial_vh)|_X$ vanish, since $X$ is reduced and contained in
$\Sing(\Theta)$.  Thus $f$ factors through $\Theta$.
Proposition~\ref{prop:theta-dual} then shows that $a_v$ factors through
$S$, so $v\in T_0S$.  As $v$ was arbitrary, $T_0S=T_0A$.
\end{proof}

\begin{proof}[Proof of Theorem~\ref{thm:main}]
A nontrivial product of principally polarized abelian varieties has a
theta divisor whose product components meet in codimension two.
Their intersection lies in its singular locus.  The hypothesis of
Theorem~\ref{thm:main} therefore forces $(B,[\Theta])$ to be
indecomposable.  Put $g=\dim B$.  For $g\leq5$ the conclusion is
Theorem~\ref{thm:beauville}, so assume
$g\geq6$.

Choose a codimension-three irreducible component
$Z\subseteq\Sing(\Theta)$.  By Theorem~\ref{thm:higher-multiplier}(i),
\begin{equation}\label{eq:hmi-containment}
 J:=I_{1,<-1/2}(\Theta)\subseteq I_Z.
\end{equation}
Choose a sufficiently small rational $0<\varepsilon<1/2$ such that
$J=I_{1,-1/2-\varepsilon}(\Theta)$.  Applying
Theorem~\ref{thm:higher-multiplier}(ii) with $k=1$ and
$\alpha=-1/2-\varepsilon$ gives
\begin{equation}\label{eq:hmi-it0}
 H^i\bigl(B,J(2\Theta)\otimes P\bigr)=0
 \quad(i>0,\ P\in\Pic^0(B)).
\end{equation}
By Theorem~\ref{thm:higher-multiplier}(iii) and
\eqref{eq:hmi-containment},
\begin{equation}\label{eq:hmi-support}
 |Z|\subseteq\Supp(\OO_B/J)\subseteq|\Sing(\Theta)|,
 \qquad \dim V(J)=g-3.
\end{equation}
In particular, $J$ is a nonzero proper coherent ideal.  Applying
Theorem~\ref{thm:bridge} and \eqref{eq:hmi-it0}, we obtain a reduced,
irreducible, geometrically nondegenerate $\GV$ subscheme
$X=V(J)$ of dimension $g-3$.  Its reducedness now upgrades the support
containment in \eqref{eq:hmi-support} to the scheme-theoretic containment
$X\subseteq\Sing(\Theta)$.

Take its theta-dual with respect to the original divisor $\Theta$.
By Theorem~\ref{thm:gv-duality}, it is an integral surface $S$.
By \eqref{eq:hmi-support} and Lemma~\ref{lem:singular-dual},
\[
 0\in S,\qquad T_0S=T_0B.
\]
Since $g\geq6>2=\dim S$, the surface is singular at $0$.
All the hypotheses of Theorem~\ref{thm:chen} are now satisfied.
Its conclusion proves Theorem~\ref{thm:main}.
\end{proof}

\section{A nonreduced base scheme in codimension three}
\label{sec:base-loci}

Theorem~\ref{thm:base} shows that every codimension-two irreducible component of the base scheme of an ample line bundle is generically reduced.
In fact, the codimension-two hypothesis is sharp: a codimension-three component
of a base scheme need not be generically reduced.  We conclude with
an explicit counterexample to show that the base scheme of an ample line bundle on an abelian variety can be non-reduced.

\begin{example}\label{ex:nonreduced-base}
In $\PP^2$ with homogeneous coordinates $X, Y, Z$, consider the
smooth plane quartic curve
\[
 C:\ X^4+Y^3Z+XZ^3=0,\qquad p=[0:1:0],\quad q=[0:0:1].
\]
Near $p$, use the affine coordinates $(t,\eta)=(X/Y,Z/Y)$;
near $q$, use $(x,y)=(X/Z,Y/Z)$.
Under these local coordinates, the  equation can be written respectively
\[
 t^4+\eta+t\eta^3=0,\qquad x^4+y^3+x=0
\]
It is clear that $t$ and $y$ are local parameters on $C$ at $p$ and $q$,
respectively.  

The line sections give $K_C\sim4p$ and $3q\sim3p$; hence
$a=[q-p]\in P:=J(C)$ has order three.  The theta divisor
$\Theta=W_2(C)-2p$ is smooth and symmetric, and contains
$0,a,-a$, represented by $2p,p+q,2q$.
Put $M=\OO_P(\Theta)$ and take the degree-three isogeny
\[
 \pi:A=(\widehat P/\phi_M(\langle a\rangle))^\vee\longrightarrow P,
 \qquad L=\pi^*M.
\]
By \cite[Lemmas~1.1--1.2]{BLR93}, $L$ has type $(1,1,3)$ and
its three eigensections form a basis of $H^0(A,L)$, with zero
divisors $\pi^{-1}(\Theta-b)$ for $b\in\langle a\rangle$.
Consequently, as schemes,
\[
 \operatorname{Bs}|L|=\pi^{-1}(W),\qquad
 W=\Theta\cap(\Theta+a)\cap(\Theta-a).
\]

We determine $W$ at $0$.  For $f(x,y)=x^4+y^3+x$, the Poincar\'e
residue description of canonical differentials
\cite[Chapter~I, Appendix~A, Exercises~8--9, p.~53]{ACGH85}
identifies $H^0(C,K_C)$ with the forms $h(x,y)\,dx/f_y$, where
$h\in\langle1,x,y\rangle$.  Since $f_y=3y^2$, we take the basis
\[
 \omega_1=\frac{dx}{y},\qquad
 \omega_2=\frac{dx}{y^2},\qquad
 \omega_3=\frac{x\,dx}{y^2}.
\]
Let $z=(z_1,z_2,z_3)$ be the additive analytic coordinates near
$0\in P$ induced by the dual basis of $\{\omega_1, \omega_2, \omega_3\}$.  Write $p(t)$ and $q(s)$ for
points near $p$ and $q$ with the indicated parameters, and the local Abel-Jacobi map is
\[
 U(t)=\left(\int_p^{p(t)}\omega_j\right)_{j=1}^3,
 \qquad
 V(y)=\left(\int_q^{q(y)}\omega_j\right)_{j=1}^3,
\]
where the integrals are taken locally. 
For example, substituting $x=t/\eta$ and $y=1/\eta$ at $p$ gives
\[
 \omega_2=\frac{dx}{y^2}
 =\eta^2\,d\!\left(\frac{t}{\eta}\right)
 =(\eta-t\eta')\,dt
 =\bigl(3t^4+O(t^{13})\bigr)\,dt.
\]
Integrating from $0$ to $t$, we obtain
$U_2(t)=\tfrac35t^5+O(t^{14})$.
Computing other terms in the same way, we have
\[
 \begin{aligned}
 U(t)&=\left(-3t+O(t^{10}),\ \tfrac35t^5+O(t^{14}),\
                         -\tfrac32t^2+O(t^{11})\right),\\
 V(y)&=\left(-\tfrac32y^2+O(y^{11}),\ -3y+O(y^{10}),\
                         \tfrac34s^4+O(s^{13})\right).
 \end{aligned}
\]

The Abel--Jacobi map $\alpha:C^{(2)}\to\Theta$, $D\mapsto[D-2p]$,
is an isomorphism since $C$ is nonhyperelliptic.
Near $2p$, write $D=p(t_1)+p(t_2)$ and use the symmetric local coordinates
$r=t_1+t_2$, $b=t_1t_2$ on $C^{(2)}$.
To express $z=U(t_1)+U(t_2)$ in these coordinates, put
$S_k=t_1^k+t_2^k$ for $k\in \mathbb N$.  The recurrence
\[
 S_0=2,\qquad S_1=r,\qquad
 S_k=rS_{k-1}-bS_{k-2}\quad(k\geq2)
\]
gives $S_2=r^2-2b$ and $S_5=r^5-5r^3b+5rb^2$.
More generally, each monomial of $S_k$ has the form $r^{k-2j}b^j$,
so $S_k\in(r,b)^{\lceil k/2\rceil}$.
Thus the above expansions give
\[
 \begin{aligned}
 z_1&=-3S_1+O(5)=-3r+O(5),\\
 z_2&=\tfrac35S_5+O(7)
      =\tfrac35\bigl(r^5-5r^3b+5rb^2\bigr)+O(7),\\
 z_3&=-\tfrac32S_2+O(6)
      =3b-\tfrac32r^2+O(6).
 \end{aligned}
\]
Here $O(k)$ denotes an element of $(r,b)^k$; below it denotes terms
of total degree at least $k$ in the coordinates under consideration.
The remainders in $U_1$, $U_2$, and $U_3$ contribute
$O(5)$, $O(7)$, and $O(6)$, respectively, and every displayed
monomial of $z_2$ has total degree at least three.  Hence $z_2\in(r,b)^3$.  The linear terms of $z_1$ and
$z_3$ show that $(z_1,z_3)$ are local coordinates on $\Theta$ at $0$,
so its equation is $z_2=F(z_1,z_3)$ with $F\in(z_1,z_3)^3$.

Near $p+q$, use $(t,y)$ for $D=p(t)+q(y)$.  Since
$\alpha(p+q)=a$, the map $\alpha-a$ near $p+q$ is $U(t)+V(s)$, so
\[
 (z_1,z_2)=(-3t,-3y)+O(2),\qquad z_3=-\tfrac32t^2+O(4).
\]
Consequently $(z_1,z_2)$ are local coordinates on $\Theta-a$ at
$0$; substituting $t=-z_1/3+O(2)$ gives its equation
$z_3=Q(z_1,z_2)$ with $Q=-z_1^2/6+O(3)$.
% Near $2q$, write $D=q(s_1)+q(s_2)$ and use
% $\rho=s_1+s_2$, $\sigma=s_1s_2$.  As $\alpha(2q)=-a$, the
% map $\alpha+a$ is locally $V(s_1)+V(s_2)$ and parametrizes
% $\Theta+a$ at $0$.  Directly taking symmetric sums gives
% \[
%  \begin{aligned}
%  z_1&=3\sigma-\tfrac32\rho^2+O(6),\qquad
%  z_2=-3\rho+O(5),\\
%  z_3&=\tfrac34\bigl(\rho^4-4\rho^2\sigma+2\sigma^2\bigr)+O(7),
%  \end{aligned}
% \]
% where the orders are taken in $(\rho,\sigma)$.
% Thus $(z_1,z_2)$ are local coordinates, and
% $z_3=z_1^2/6+O(3)$ after eliminating $\rho,\sigma$.
Since $-(\Theta-a)=\Theta+a$, we therefore obtain
\[
 \begin{aligned}
 \Theta&:\quad z_2=F(z_1,z_3),\qquad F\in(z_1,z_3)^3,\\
 \Theta-a&:\quad z_3=Q(z_1,z_2),\qquad
 Q=-\tfrac16z_1^2+O(3),\\
 \Theta+a&:\quad z_3=-Q(-z_1,-z_2).
 \end{aligned}
\]
The Jacobian of the first two equations with respect to $(z_2,z_3)$
is the identity at $0$.  They therefore define a smooth formal curve
with parameter $u=z_1$, on which $z_2=O(u^3)$ and
$z_3=-u^2/6+O(u^3)$.  On this curve, the remaining equation is
\[
 z_3+Q(-z_1,-z_2)=-\tfrac13u^2+O(u^3)=0.
\]
Thus
\[
 \widehat{\OO}_{W,0}\simeq\CC[[u]]/(u^2).
\]
Translation by $a$ preserves $W$, so $0,a,-a$ are isolated double
points of $W$ and in fact, they exhaust $W$.
% To see that they exhaust $W$, let
% $D=\Theta\cap(\Theta-a)$.  The theta divisors are distinct and
% irreducible, so $D$ is a pure complete-intersection curve.
% Write $[D]=\sum_i m_i[D_i]$ for its fundamental cycle.
% The three isolated intersections with $\Theta+a$ have total length
% six; by associativity of intersection multiplicities,
% \[
%  6\leq\sum_{D_i\not\subset\Theta+a}m_i\deg(M|_{D_i})
%  \leq M\cdot[D]=\Theta^3=6.
% \]
% Ampleness makes every $\deg(M|_{D_i})$ positive, so no component
% of $D$ is contained in $\Theta+a$.  Thus $W$ is a zero-dimensional
% complete intersection of length six, already exhausted by these three
% double points.  
Since $\pi$ is \'etale of degree three,
\[
 \operatorname{Bs}|L|\simeq
 \coprod_{j=1}^{9}\operatorname{Spec}\CC[\varepsilon]/(\varepsilon^2).
\]
Finally, for a globally generated ample line bundle $N$ on an
abelian variety $E$, the equality
\[
 \operatorname{Bs}|L\boxtimes N|=\operatorname{Bs}|L|\times E
\]
produces generically nonreduced components of codimension three
in every dimension at least three.

Rubei \cite[Remark~3, p.~367]{Rub98} observed that a zero-dimensional
base locus of type $(1,1,3)$ has either nine or eighteen points.
The present example realizes the nine-point case, with its nonreduced
scheme structure determined by the calculation above.
Related constructions from plane quartics and points of order three
appear in \cite[Theorem~1(1c)]{Rub99}.
\end{example}


\begin{thebibliography}{99}

\bibitem{AM67}
A.~Andreotti and A.~L.~Mayer,
\emph{On period relations for abelian integrals on algebraic curves},
Ann. Scuola Norm. Sup. Pisa (3) \textbf{21} (1967), 189--238.
\href{https://numdam.org/item/ASNSP_1967_3_21_2_189_0/}{Numdam}.

\bibitem{ACP21}
E.~Arbarello, G.~Codogni, and G.~Pareschi,
\emph{Characterizing Jacobians via the KP equation and via flexes and
degenerate trisecants to the Kummer variety: an algebro-geometric approach},
J. Reine Angew. Math. \textbf{777} (2021), 251--271.
\href{https://doi.org/10.1515/crelle-2021-0020}{\nolinkurl{doi:10.1515/crelle-2021-0020}}.
See the withdrawal notice \cite{ACP25}.

\bibitem{ACP25}
E.~Arbarello, G.~Codogni, and G.~Pareschi,
\emph{Characterizing Jacobians via the KP
equation and via flexes and degenerate trisecants to the Kummer variety:
an algebro-geometric approach}, withdrawn version and withdrawal notice,
\href{https://arxiv.org/abs/2009.14324v4}{arXiv:2009.14324v4},
3 June 2025.

\bibitem{ACGH85}
E.~Arbarello, M.~Cornalba, P.~A.~Griffiths, and J.~Harris,
\emph{Geometry of Algebraic Curves}, Volume~I,
Grundlehren der mathematischen Wissenschaften, vol.~267,
Springer, New York, 1985.
\href{https://doi.org/10.1007/978-1-4757-5323-3}{\nolinkurl{doi:10.1007/978-1-4757-5323-3}}.

\bibitem{Beauville77}
A.~Beauville,
\emph{Prym varieties and the Schottky problem},
Invent. Math. \textbf{41} (1977), 149--196.
\href{https://doi.org/10.1007/BF01418373}{\nolinkurl{doi:10.1007/BF01418373}}.

\bibitem{BL04}
C.~Birkenhake and H.~Lange,
\emph{Complex Abelian Varieties}, second ed.,
Grundlehren der mathematischen Wissenschaften, vol.~302,
Springer, Berlin, 2004.

\bibitem{BLR93}
C.~Birkenhake, H.~Lange, and S.~Ramanan,
\emph{Primitive line bundles on abelian threefolds},
Manuscripta Math. \textbf{81} (1993), 299--310.
\href{https://doi.org/10.1007/BF02567860}{\nolinkurl{doi:10.1007/BF02567860}}.

\bibitem{CMPS18}
S.~Casalaina-Martin, M.~Popa, and S.~Schreieder,
\emph{Generic vanishing and minimal cohomology classes on abelian fivefolds},
J. Algebraic Geom. \textbf{27} (2018), no.~3, 553--581.
\href{https://people.math.harvard.edu/~mpopa/papers/5GV.pdf}{Authors' manuscript}.

\bibitem{Chen26}
Y.~Chen,
\emph{Generic vanishing subschemes of codimension three},
\href{https://arxiv.org/abs/2609.15391v2}{arXiv:2609.15391v2},
17 September 2026.

\bibitem{Debarre88}
O.~Debarre,
\emph{Sur les vari\'et\'es ab\'eliennes dont le diviseur th\^eta est
singulier en codimension~3},
Duke Math. J. \textbf{57} (1988), no.~1, 221--273.
\href{https://perso.imj-prg.fr/olivier-debarre/wp-content/uploads/debarre-pub/7.pdf}{Author's copy}.

\bibitem{Debarre95}
O.~Debarre,
\emph{Minimal cohomology classes and Jacobians},
J. Algebraic Geom. \textbf{4} (1995), 321--335.
\href{https://perso.imj-prg.fr/olivier-debarre/wp-content/uploads/debarre-pub/28.pdf}{Author's copy}.

\bibitem{Debarre95Connectivity}
O.~Debarre,
\emph{Th\'eor\`emes de connexit\'e et vari\'et\'es ab\'eliennes},
Amer. J. Math. \textbf{117} (1995), no.~3, 787--805.
\href{https://doi.org/10.2307/2375089}{\nolinkurl{doi:10.2307/2375089}}.

\bibitem{Debarre99}
O.~Debarre,
\emph{Tores et vari\'et\'es ab\'eliennes complexes},
Cours Sp\'ecialis\'es, vol.~6,
Soci\'et\'e Math\'ematique de France, Paris, and EDP Sciences, Les Ulis, 1999.
\href{https://perso.imj-prg.fr/olivier-debarre/wp-content/uploads/debarre-pub/csp18.pdf}{Author's revised version}.
Page references are to the published edition.

\bibitem{Debarre06}
O.~Debarre,
\emph{Singularities of divisors on abelian varieties},
lecture notes, 20 March 2006.
\href{https://www.math.ens.psl.eu/~debarre/DivVarAb.pdf}{Original notes};
\href{https://perso.imj-prg.fr/olivier-debarre/wp-content/uploads/sites/34/2025/05/DivVarAb2024.pdf}{author's updated version}.

\bibitem{EL97}
L.~Ein and R.~Lazarsfeld,
\emph{Singularities of theta divisors and the birational geometry of
irregular varieties},
J. Amer. Math. Soc. \textbf{10} (1997), 243--258.

\bibitem{GL11}
M.~G.~Gulbrandsen and M.~Lahoz,
\emph{Finite subschemes of abelian varieties and the Schottky problem},
Ann. Inst. Fourier (Grenoble) \textbf{61} (2011), no.~5, 2039--2064.
\href{https://doi.org/10.5802/aif.2665}{\nolinkurl{doi:10.5802/aif.2665}}.

\bibitem{Horing07}
A.~H\"oring,
\emph{$M$-regularity of the Fano surface},
C. R. Math. Acad. Sci. Paris \textbf{344} (2007), no.~11, 691--696.
\href{https://doi.org/10.1016/j.crma.2007.04.008}{\nolinkurl{doi:10.1016/j.crma.2007.04.008}}.

\bibitem{Horing10}
A.~H\"oring,
\emph{Minimal classes on the intermediate Jacobian of a generic cubic threefold},
Comm. Contemp. Math. \textbf{12} (2010), no.~1, 55--70.
\href{https://doi.org/10.1142/S0219199710003737}{\nolinkurl{doi:10.1142/S0219199710003737}}.

\bibitem{JP20}
Z.~Jiang and G.~Pareschi,
\emph{Cohomological rank functions on abelian varieties},
Ann. Sci. \'Ec. Norm. Sup\'er. (4) \textbf{53} (2020), 815--846.
\href{https://doi.org/10.24033/asens.2435}{\nolinkurl{doi:10.24033/asens.2435}}.

\bibitem{LT15}
L.~Lombardi and S.~Tirabassi,
\emph{$GV$-subschemes and their embeddings in principally polarized
abelian varieties},
Math. Nachr. \textbf{288} (2015), nos.~11--12, 1405--1412.
\href{https://doi.org/10.1002/mana.201400238}{\nolinkurl{doi:10.1002/mana.201400238}}.

\bibitem{Meng26}
F.~Meng,
\emph{Generic vanishing conjecture for divisors},
\href{https://arxiv.org/abs/2609.11377v1}{arXiv:2609.11377v1},
10 September 2026.

\bibitem{Mukai81}
S.~Mukai,
\emph{Duality between $D(X)$ and $D(\widehat X)$ with its application
to Picard sheaves},
Nagoya Math. J. \textbf{81} (1981), 153--175.

\bibitem{PP08}
G.~Pareschi and M.~Popa,
\emph{Generic vanishing and minimal cohomology classes on abelian varieties},
Math. Ann. \textbf{340} (2008), 209--222.
\href{https://doi.org/10.1007/s00208-007-0146-7}{\nolinkurl{doi:10.1007/s00208-007-0146-7}}.

\bibitem{Ran81}
Z.~Ran,
\emph{On subvarieties of abelian varieties},
Invent. Math. \textbf{62} (1981), 459--479.

\bibitem{Rub98}
E.~Rubei,
\emph{Projective normality of abelian varieties with a line bundle of
 type $(2,\ldots)$},
Boll. Unione Mat. Ital. B (8) \textbf{1} (1998), no.~2, 361--367.
\href{https://www.bdim.eu/item?id=BUMI_1998_8_1B_2_361_0}{BDIM}.

\bibitem{Rub99}
E.~Rubei,
\emph{On normal generation of line bundles on abelian threefolds},
Kodai Math. J. \textbf{22} (1999), no.~1, 15--34.
\href{https://doi.org/10.2996/kmj/1138043985}{\nolinkurl{doi:10.2996/kmj/1138043985}}.

\bibitem{Schreieder16}
S.~Schreieder,
\emph{Theta divisors with curve summands and the Schottky problem},
Math. Ann. \textbf{365} (2016), nos.~3--4, 1017--1039.
\href{https://doi.org/10.1007/s00208-015-1287-8}{\nolinkurl{doi:10.1007/s00208-015-1287-8}}.

\bibitem{Schreieder17}
S.~Schreieder,
\emph{Decomposable theta divisors and generic vanishing},
Int. Math. Res. Not. IMRN \textbf{2017} (2017), no.~16, 4984--5009.
\href{https://doi.org/10.1093/imrn/rnw160}{\nolinkurl{doi:10.1093/imrn/rnw160}}.

\bibitem{SY26}
C.~Schnell and R.~Yang,
\emph{Higher multiplier ideals},
J. Reine Angew. Math. \textbf{832} (2026), 153--197.
\href{https://doi.org/10.1515/crelle-2025-0097}{\nolinkurl{doi:10.1515/crelle-2025-0097}}.
See also \href{https://arxiv.org/abs/2309.16763v7}{arXiv:2309.16763v7}.

\end{thebibliography}
\end{document}